\documentclass [12pt]{amsart}
\usepackage{amsmath}%
\usepackage{amsfonts}
\usepackage{amssymb}%
\usepackage{graphicx}

\newtheorem{theorem}{Theorem}[section]

\newtheorem{corollary}[theorem]{Corollary}

\newtheorem{definition}[theorem]{Definition}

\newtheorem{lemma}[theorem]{Lemma}

\newtheorem{remark}[theorem]{Remark}

\newtheorem*{thm1*}{Theorem 1}
\begin{document}

\title{$Q_{\alpha}$-Normal Families and entire functions}
\author [Shai Gul and Shahar Nevo]{Shai Gul and Shahar Nevo }

\begin{abstract}
For every countable ordinal number $\alpha$ we construct an entire function $f=f_\alpha$ such that the family $ \left\{f(nz):n\in \mathbb {N}\right\}$ is exactly $Q_\alpha$-normal in the unit disk.
\end{abstract}
\thanks{The authors recieved support from the Israel Science Foundation, Grant No. 395/07. This research is part of the European Science Networking Programme HCAA}

\maketitle \tableofcontents
\section {Introduction}
\setcounter{equation}{0} \numberwithin{equation}{section}

The notion of $Q_m$-normal family of meromorphic functions
$(m\in N)$ was developed by C.T Chaung \cite{Chuang}. This is a geometric extention to the well-known notions of normality and  quasi-normality due to P. Montel. A family $F$ of functions, meromorphic on a domain $D\subset C$, is called $Q_m$-normal in $D$, if from each sequence $S \subset F$ we can extract a subsequence $\hat S $, such that $\hat S $ converges uniformaly with respect to the spherical metric $ \chi $ to limit function $ f $ on a domain $ D\backslash E $, where $ E \subset D$ satisfies $ E_D^{(m)}  = \emptyset $.
Here $f$ is meromorphic or $f\equiv\infty$ and  $E_D^{(m)}$ is the derived set of order $m$ of $E$ with respect to $D$. When $m=0$ then $ F $ is a normal family and when $ m=1 $ then $F$ is a quasi-normal family.
Results in this subject were achieved in \cite{Applications},\cite{Generating} and \cite{Transfinite}, and some of them will be detailed in the sequel. In particular, in \cite{Transfinite}, this notation of $Q_m$-normality was extended further to $Q_\alpha$-normality, where $\alpha$ is any ordinal number, in an analogous way that will be explained.

The goal of this paper is to continue the progress of this research and to add foundations for the continuing study of $Q_\alpha$-normal families. We will derive a few elementary results that are summarized in the introductory section, but our main result is an extension of a result from  \cite{Generating}. For an abritrary countable ordinal number $\alpha$ we shall construct an entire function $f(z)=f_\alpha(z)$, such that the family $F(f):=\left\{{f(nz):n\in \mathbb{N}}\right\}$  is exactly $Q_\alpha$-normal in the unit disc $\Delta$.
For ordinal numbers having an immediate predecessor (i.p), the result is sharper. Because of the inductive nature of the definition of a $Q_\alpha$-normal family, transfinite induction plays a major role throughout this text.

This paper is organized as follows. In Section 2 we bring background on $Q_\alpha$-normality, and usually we do not give proofs. However, the new basic results will be proved and also a few of the old ones, in cases where we want to emphasis the difference between dealing with ordinal numbers having i.p and those having limit ordinal numbers.

We shall also state our main result. In Section 3 we prove the main theorem for the case of $Q_\alpha$-normal families of finite order, where $\alpha$ has an i.p. In Section 4 we prove the limit cases of the theorem, i.e. a $Q_\alpha$-normal family of infinite order where $\alpha$ has an i.p or the case where $\alpha$ is a limit ordinal number (without relating to the order). The notion of ``order" of normal family is explained in Section 2. In Section 5 we make a few remarks about the functions $f_\alpha(z)$.

For an arbitrary ordinal number $\alpha$ we use the usual notation $ \Gamma (\alpha ): = \left\{ {\beta :\beta  < \alpha } \right\} $.
We shall use the following notation and conventions.\newline
\newline
1.1) Let $z_0  \in \mathbb{C}, r > 0$.
\begin{align*}
\begin{array}{l}
 \Delta \left( {z_0 ,r} \right) = \left\{ {\left| {z - z_0 } \right| < r} \right\}. \\
 \bar \Delta \left( {z_0 ,r} \right) = \left\{ {\left| {z - z_0 } \right| \le r} \right\}. \\
 \Delta '\left( {z_0 ,r} \right) = \left\{ {0 < \left| {z - z_0 } \right| < r} \right\}. \\
 \end{array}
 \end{align*}
\newline
1.2) Let $A\subset \mathbb{C}$, $z\in \mathbb{C} $. Then
$z \cdot A = \left\{ {z \cdot a:a \in A} \right\}$.\newline
\newline
1.3) Let $D\subset \mathbb{C}$ be a domain and $A\subset D$. We say that $A$ is compactly contained in $D$ if $\overline{A} \subset D$.\newline
\newline
1.4) Let $A_1,A_2,...$ be sets in $\mathbb{C}$. These sets are said to be strongly disjoint if $\overline{A_i}\cap\overline{A_j}=\emptyset$ whenever $i\neq j$.\newline
\newline
1.5) Let $\left\{g_n\right\}$ be a sequence of meromorphic functions on a domain $D\subset \mathbb{C}$. If $\left\{g_n\right\}$ converges uniformly on a compact subset of $D$ to $g$ with respect to the spherical metric $\chi$ on $\hat{\mathbb{C}}$ (then $g$ is a meromorphic function on $D$ or the constant $\infty$), we say that $\left\{g_n\right\}$ converges to $g$ locally $\chi$-uniformaly on $D$ and write $g_n \mathop  \Rightarrow \limits^\chi  g $ on $D$.

In case the functions $g_n$ are holomorphic in $D$, then either the convergence  is locally uniform with respect to the Eucliden metric, in which case the limit function $g$ is holomorphic on $D$, or $\left\{g_n\right\}$ diverges uniformly to $\infty$ on compacta. In this case, we write $g_n\Rightarrow g$ on $D$ or $g_n\Rightarrow \infty$ on $D$, respectively.\newline
\newline
1.6) Let $0\leq\theta<2\pi$, $0<\varepsilon<\pi$. The infinite open angular sector of the opening $2\varepsilon$ around $\theta$ is $S(\theta,\varepsilon)=\left\{z:\theta-\varepsilon<arg{z}<\theta+\varepsilon\right\}$. The ray $L(\theta)$ is $\left\{z:arg{z}=\theta \right\}$.

In Section 6 we introduce an open problem concerning the first uncountable ordinal number.
\section {Background on $Q_\alpha$-normal families of meromorphic functions}
Many of the definitions and statements in this section have analogues in \cite{Chuang} and \cite{Applications} for $Q_m$-normal families or appears also in \cite{Transfinite}.
\subsection {Geometrical Background with connection to ordinal numbers}
\begin{definition}
(\cite[Definition 2.8]{Transfinite}, cf. \cite[p.163]{Kura.})
Let $D$ be a domain in the complex plane and $E\subset D$ a point set. The derived set of the first order (order $1$) of $E$ with respect to $D$, denoted by $E_D^{(1)}$, is the set of all accumulation points of $E$ in $D$, that is, $E_D^{(1)}  = E' \cap D$. For every $\alpha\geq 2$ which has an i.p, define the derived set of order $\alpha$ of $E$ with respect to $D$ by
$E_D^{(\alpha )}  = (E_D^{(\alpha  - 1)} )_D^{(1)}$. If $\alpha$  is a limit ordinal define
$E_D^{(\alpha )}  = \mathop  \bigcap \limits_{1 \le \beta  < \alpha } E_D^{(\beta )}$. We also set $E_D^{(0)}=E$. More generally, we write $E_A^{(1)}  = E' \cap A$ where $A$ is an arbitrary subset of $\mathbb{C}$ and define $E_A^{(\alpha)}$ in a similar fashion.

\end{definition}
The following is easy to prove.
\begin{corollary}  \cite [Corollary 2.9]{Transfinite} 
If $\bar E \subset A \subset \mathbb{C}$ and $\bar E \subset B \subset \mathbb{C}$, then $E_A^{(\alpha)}=E_B^{(\alpha)}$ for every ordinal number $\alpha$.
\end {corollary}
\begin{lemma}  \cite [Lemma 2.10]{Transfinite} 
Let $E\subset D$ be a point set. Then $E_D^{(\beta)}\subset E_D^{(\alpha)}$ for any two ordinals numbers $\alpha$, $\beta$ such that $\beta\geq \alpha \geq 1$.
\end {lemma}
This lemma is well-known, so we omit the proof.

\begin {lemma}  \cite[Lemma 2.12]{Transfinite} 
Let $\alpha_1,\alpha_2,...$ be a sequence of countable ordinal numbers. Then there is a countable ordinal number $\alpha_0$ such that for each $k\geq1$, $\alpha_k<\alpha_0$.
\end {lemma}
This follows easily from the uncountability of the set $\Gamma(\Omega)$, where $\Omega$ is the first uncountable ordinal number.
\begin{lemma}  \cite[Lemma 2.13]{Transfinite} 
Let $E\subset D$ a point set. Then there is a countable ordinal number $\alpha$ such that $E_D^{(\beta)}=E_D^{(\alpha)}$ for every $\beta\geq\alpha$.
\end{lemma}
This Lemma is Ex.7 on p.163 of \cite{Kura.}, see also \cite[Lemma 2.13]{Transfinite}.
\begin{lemma} 
If $E\subset D$ and $E=E_D^{(1)}$ then $E=\emptyset$ or $\left|E\right|= \aleph$.
\end{lemma}
This lemma is well-known, and in case of the second possibility, $E$ is called a perfect set.

\begin{lemma} 
Let $E \subset D$ and assume that $E_D^{(\alpha)}= \emptyset $ for some countable ordinal number $\alpha$. Then $E$ is countable.
\end{lemma}
\begin{proof}
The lemma holds for $\alpha=0$ or $\alpha=1$. We apply now transfinite induction. Assume that the lemma is true for every $\beta<\alpha$. If $\alpha$ has an i.p, $\alpha-1$, then $E_D^{(\alpha-1)}$ is a discrete set in $D$ (and of course countable). Thus
$\hat D: = D{\rm{\backslash E}}_D^{(\alpha  - 1)} $ is a domain, and then we set $\hat E = E \backslash {\rm{E}}_D^{(\alpha  - 1)}$.
Clearly $\hat E_{\hat D}^{\left( {\alpha  - 1} \right)}  = \emptyset $ and by the induction assumption $\hat E $ is countable. Since $E\subset\hat{E}\cup E_D^{(\alpha-1)}$, then also $E$ is countable. Suppose now that $\alpha$ is a limit ordinal number and let $z\in D$. By Definition 2.1 there is some $\beta<\alpha$ such that $z\notin E_D^{(\beta+1)}$. Hence there is some $r_z >0 $ such that $\Delta '(z_0 ,r) \cap E_D^{(\beta )}  = \emptyset$,
so if we denote $E_z  = E \cap \Delta (z,r)$, then
$(E_z)_{\Delta (z_0 ,r)} ^{(\beta  + 1)}  = \emptyset$
and by the the induction assumption $E_z=\emptyset$.
Now, since every collection of open sets has a countable sub-collecion having the same union, we deduce that $E$ is countable.
\end{proof}
From this lemma and from Lemma 2.5, we deduce
\begin{corollary}   
If $E\subset D$ is not counatable, then $ \left| {\mathop  \bigcap \limits_\alpha ^{} E_D^{(\alpha )} } \right| = \aleph $ (cf. \cite[Lemma 2.13]{Applications}).
\end{corollary}
The following lemma is obvious.
\begin{lemma} 
If $E \subset D$ and $ \alpha $ is an ordinal number and $z_0 \in E_D^{(\alpha)}$, then $z_0 \in E_A^{(\alpha)}$ for any open set $ A \subset D$ such that $z_0 \in A$.
\end {lemma}
\begin {definition} (\cite[Definition 9.3]{Chuang},\cite[Definition 2.14]{Transfinite}) 
Let $E_n$ $(n=1,2,\dots)$ and $E$ be sets of points in the complex plane $\mathbb{C}$. We say that $E$ is a limit set of the sequence $E_n$ $(n=1,2,\dots)$ if for any point $z_0 \in E $, any positive number $\varepsilon$, and any positive integer $N$, we can find an integer $n \geq N $ such that
$E_n  \cap \Delta \left( {z_0 ,\varepsilon } \right) \ne \emptyset$.
\end {definition}
\begin{lemma}  \cite [Lemma 2.15]{Transfinite} 
Let  $\left\{ {A_n } \right\}_{n = 1}^\infty$ be a sequence of sets belonging to a domain $D$. Then ${A_n}$ has a non empty limit set $E$, if and only if there exists a sequence of points $\left\{ {a_k } \right\}_{k = 1}^\infty $, $a_k  \in A_{n_k } $ ($n_1<n_2<..$) such that
$a_k \mathop  \to \limits_{k \to \infty } a_0$ for some $a_0 \in D $.
\end {lemma}
The proof of this lemma is obvious.
\begin{lemma} \cite[Theorem 2.16]{Transfinite} 
Let $\left\{ {A_n } \right\}_{n = 1}^\infty  $ be a sequence of subsets of $D$ whose limit set in $D$ is empty. Then for each ordinal $\alpha$,
\begin{equation}
\left( {\bigcup\limits_{n = 1}^\infty  {A_n } } \right)_D^{(\alpha )}  = \bigcup\limits_{n = 1}^\infty  {(A_n )_D ^{\left( \alpha  \right)} }.
\end{equation}
\end{lemma}
\begin{proof}
By Lemma 2.4 and Lemma 2.5, it is enough to prove the lemma in the case where $\alpha$ is countable. It is clear that the right side is contained in the left side. We have to show the opposite containment. We proceed by using transfinite induction. Suppose that $\alpha=1$ and that
$z_0  \in \left( {\bigcup\limits_{n = 1}^\infty  {A_n } } \right)_D^{(1)}$.
If
$z_0  \notin \left( {\bigcup\limits_{n = 1}^\infty  {A_n } } \right)_D^{(1)}$
we deduce the existence of sequence of points $\left\{ {a_k } \right\}_{k = 1}^\infty  $, $n_1 < n_2 < \dots{}$; $a_k \in A_{n_k} $ for each $k \geq 1$ such that $
a_k \mathop  \to \limits_{k \to \infty } z_0 $. By Lemma 2.11, $A_n$ has a nonempty limit set, a contradiction. Here we have proved the theorem for the case $\alpha=1$.
Assume that the theorem is true for $ \beta<\alpha $, where $\alpha$ is an ordinal number. As usual we consider two cases:\newline
\newline
\textbf{(1) $\alpha-1$ exists}. By the induction assumption,
\begin{equation}
\left( {\mathop  \bigcup \limits_{n = 1}^\infty  A_n } \right)_D^{(\alpha )}  = \left( {\mathop  \bigcup \limits_{n = 1}^\infty  (A_n )_D^{(\alpha  - 1)} } \right)_D^{(1)}.
\end{equation}
We claim that
$\left\{ {(A_n )_D^{(\alpha  - 1)} } \right\}_{n = 0}^\infty $
has no nonempty limit set in $D$. Indeed, if this is not the case, then by Lemma 2.11 we would have $z_k\rightarrow z_0$,
$z_k  \in (A_n )_D^{(\alpha  - 1)}$. Then by Lemma 2.3, $z_k  \in (A_n )_D^{(1 )}$, $k=1,2,...$, and we could construct a sequence $ \left\{ {\zeta _k } \right\}_{k = 1}^\infty $, $ \zeta_k \in A_{n_k} $ such that $\zeta_k \rightarrow z_0$. This contradicts the assumption that ${A_n}$ has no nonempty limit set. Now, by the case $\alpha=1$, the right side of (2.2) is $\mathop  \bigcup \limits_{n = 1}^\infty  (A_n )_D^{(\alpha )}$, as desired.\newline
\newline
\textbf{(2) $\alpha$ is a limit ordinal number}. We have by Definition 2.1 and the induction assumption, $\left( {\mathop  \bigcup \limits_{n = 1}^\infty  A_n } \right)_D^{(\alpha )}  = \mathop  \bigcap \limits_{\beta  < \alpha } \left( {\mathop  \bigcup \limits_{n = 1}^\infty  A_n } \right)_D^{(\beta )}  = \mathop  \bigcap \limits_{\beta  < \alpha } \mathop  \bigcup \limits_{n = 1}^\infty  (A_n )_D^{(\beta )}$; and for the right side of (2.1), we have
$\bigcup\limits_{n = 1}^\infty  {\left( {A_n } \right)_D^{\left( \alpha  \right)} }  = \bigcup\limits_{n = 1}^\infty  {\bigcap\limits_{\beta  < \alpha } {\left( {A_n } \right)_D^{\left( \beta  \right)} } }$.
Thus we need to show that

\begin{equation}
\bigcap\limits_{\beta  < \alpha } {\bigcup\limits_{n = 1}^\infty  {\left( {A_n } \right)} _D^{\left( \beta  \right)}  \subset \bigcup\limits_{n = 1}^\infty  {\bigcap\limits_{\beta  < \alpha } {\left( {A_n } \right)_D^{\left( \beta  \right)} } } }.
\end{equation}
Let $z_0$ belong to the left side of (2.3). For every $ \beta<\alpha $,
$z_0  \in \left( {A_n } \right)_D^{(\beta )}$
for a finite number of values of $n$ (otherwise in a manner similar to the previous case, it can be shown that $z_0$ is a limit point of $\left\{ {A_n } \right\}$, a contradiction).
In particular, there is some $ n_0\in \mathbb{N}$ such that
$z_0  \notin \left( {A_n } \right)_D^{(1)}$
for $n > n_0$, and by Lemma 2.3
\begin{equation}
\begin{array}{l}
 z_0  \in \mathop  \bigcap \limits_{\beta  < \alpha } \mathop  \bigcup \limits_{n = 1}^n \left( {A_n } \right)_D^{(\beta )}  \\

 \end{array}
 \end{equation}
Denote for $1 \le n{\rm{ }} \le n_0  $,
$B_n  = \{ \beta  < \alpha :z_0  \in (A_n )_D^{(\beta )} \}$. Then by (2.4) $\Gamma (a) = \mathop  \cup \limits_{n = 1}^{n_0 } B_n$.
Thus, for some $k$, $1 \le k \le n_0$, $B_k$ contains ordinal numbers ``as large as we want", and by Lemma 2.3
$z_0  \in \bigcap\limits_{\beta  < \alpha }^{} {\left( {A_k } \right)_D^{(\beta )} }$
and $z_0$ is in the right side of (2.3). The proof is completed.
\end{proof} 
As a corollary of Lemma 2.12 we have
\begin {corollary}  \cite [Corrollary 2.17]{Transfinite} 
If $A,B \subset D$ and $\alpha$ is an ordinal number, then
$\left( {A \cup B} \right)_D^{\left( \alpha  \right)}  = A_D^{\left( \alpha  \right)}  \cup B_D^{\left( \alpha  \right)} $.
\end {corollary}

\subsection{Background on $Q_\alpha$-normal families}
Using the geometrical background we proceed now to $Q_\alpha$ normal families.
\begin{definition}  (\cite [Definition 2.1]{Transfinite}, cf. \cite[Definition 1.4]{Chuang},\cite [Definition 3.1]{Applications}) 
Let $S=\left\{{f_n}\right\}$ be a sequence of meromorphic function in a domain $D$. A point $z_0$ of $D$ is called a $C_0$-point of $S$, if there is a disk $\Delta(z_0,r)$ contained in $D$ such that the sequence $S$ is uniformly convergence in $\Delta(z_0,r)$ with respect to the spherical distance. The sequence $S$ is said to be a $C_0$-sequence in $D$ if each point in $D$ is a $C_0$-point of $S$.
\end{definition}
If $z_0$ is not a $C_0$-point of $S$ then $z_0$ is called a non$C_0$-point of $S$.
Observe that the set $E\subset D$ of non$C_0$-points of $S$ is a closed set relative to $D$.\newline
\newline
We can now define a $C_\alpha$-point for any ordinal number $\alpha$.
\begin{definition} (\cite[Definitions 2.2, 2.5]{Transfinite}, cf. \cite[Definitions 8.1, 8.2]{Chuang},\cite[Definitions 3.15, 3.17]{Applications}) 
Let $\alpha\geq1$ be an ordinal number. Let $S=\left\{{f_n}\right\}$ be a sequence of meromorphic functions in a domain $D$ and $z_0$ a point of $D$. We say thet $z_0$ is a $C_\alpha$-point of $S$ if the following conditions hold:\newline
(1) If $\alpha$ has an i.p $\alpha-1$, then there exists $r>0$ such that $ \Delta(z_0,r)\subset D $ and each point in $\Delta'(z_0,r)$ is a $C_{\alpha-1}$-point of $S$.\newline
(2) If $\alpha$ is a limit ordinal, then $z_0$ is a $C_\beta$-point of $S$ for some $ \beta<\alpha$.

If $z_0$ is not a $C_\alpha $-point of $S$, then $z_0$ is called a non$C_\alpha$-point of $S$.
By transfinite induction, it is readily seen the the notion of $C_\alpha$-point is well defined for every ordinal $\alpha$.
\end{definition}
A sequence $S$ of meomorphic functions on $D$ is called a $C_\alpha$-sequence, if each point of $D$ is a $C_\alpha$-point of $S$.
\begin{lemma} (cf. \cite[Lemmas 2.4, 2.6]{Transfinite} ,\cite [Lemmas 8.1, 8.2]{Chuang}, \cite [Lemmas 3.16, 3.17]{Applications}) 
Suppose that $z_0\in D $ is a $C_\alpha$-point of a sequence $S=\left\{{f_n}\right\}$. Then $z_0$ is a $C_\beta$-point of $S$ for every $\beta>\alpha$.
\end{lemma}
Thus if $S$ is a $C_\alpha$-sequence, then $S$ is also a $C_\beta$-sequence for every $\beta>\alpha$.
\begin{definition} (\cite[Definition 2.21]{Transfinite}, cf. \cite[Definition.8.5]{Chuang}, \cite[Definition 3.24]{Applications}) 
Let $F$ be a family of meoromorphic functions in a domain $D$ and $\alpha$ an ordinal number. We say that the family $F$ is $Q_\alpha$-normal in $D$, if from every sequence of functions of the family $F$, we can extract a subsequence which is a $C_\alpha$-sequence in $D$. $F$ is said to be $Q_\alpha$-normal at a point $z_0$ of $D$ if there is a disk $\Delta(z_0,r)$ contained in $D$ such that $F$ is $Q_\alpha$-normal in $\Delta(z_0,r)$. (In particular, a $Q_0$-normal family is a normal family and a $Q_1$-normal family is a quasi-normal family, respectively.)
\end{definition}
Evidently, according to Definition 2.17 if $F$ is $Q_\alpha$-normal in $D$, then $F$ is $Q_\alpha$-normal at each point $D$. Conversley, we have the following theorem.
\begin{theorem}(\cite[Theorem 2.22]{Transfinite}, cf. \cite [Theorem 8.3]{Chuang}, \cite [Theorem 3.25]{Applications}) 
Let $F$ be a family of meromorphic functions in a domain $D$ and $\alpha$ be an ordinal number. If $F$ is $Q_\alpha$-normal at each point of $D$, then $F$ is $Q_\alpha$-normal in $D$.
\end{theorem}
The proof goes exactly as the proof of Theorem 1.4 in \cite{Chuang}. One simply writes everywhere ``$Q_\alpha$-normal" instead of ``normal" and ``$C_\alpha$ -sequence",  instead of ``$C_0$-sequence".

A trivial consequence of Lemma 2.15 is
\begin{lemma}  \cite [Lemma 2.23]{Transfinite} 
If for some ordinal number $\alpha$, a family $F$ of meomorphic functions in a domain $D$ is $Q_\alpha$-normal in $D$, then $F$ is $Q_\beta$-normal in $D$ for every $\beta>\alpha$.
\end{lemma}
We give now a few useful properties of $C_\alpha$-sequences.
\begin{lemma} (\cite[Lemma 2.20]{Transfinite}, cf. \cite[Remark 1.27]{Applications}) 
Let $S=\left\{{f_n}\right\}$ be a sequence of meromorphic functions in a domain $D$, $z_0\in D$, and let $\alpha$ be an ordinal number. Then $z_0$ is a non$C_\alpha$-point of $S$ if and only if $z_0\in E_D^{(\alpha)}$, where $E$ is the set of non$C_0$-points of $S$ in $D$.
\end{lemma}
The proof of this lemma follows easily by transfinite induction. As a consequence we have
\begin{lemma} (\cite[Theorem 2.19]{Transfinite}, cf.\cite[Theorem 8.2]{Chuang}) 
Let $S$ be a sequence of meromorhic functions in a domain $D$, and $\alpha$ an ordinal number. In order for $S$ to be a $ C_\alpha$-sequence in $D$, it is necessary and sufficient that the set of non$C_0$-points of $S$ in $D$ satisfy $E_D^{(\alpha)}=\emptyset$.
\end{lemma}
By Lemma 2.20 and Lemma 2.7 we deduce that the set of non$C_0$-points of some $C_\alpha$-sequence, when $\alpha$ is countable ordinal number, is at most enumerable. The converse is also true.
\begin{lemma}
Let $S$ be a sequence of meromorphic functions in $D$, and let $E$ be the set of non$C_0$-points of $S$ in $D$. Then $E$ is at most countable if and only if $S$ is a $C_\alpha$-sequence for some countable ordinal number $\alpha$.
\end{lemma}
\begin{proof}
We only have to prove sufficiency. Assume that $E$ is at most countable. Since $E$ is closed in $D$, then by Lemma 2.3 $E_D^{(\beta)}\subset E$ for every ordinal number $\beta$. By Lemma 2.5 there is some countable ordinal number $ \alpha $ such that $ E_D^{(\alpha )}  = E_D^{(\beta )} $ for every $\beta \geq \alpha$. If $ E_D^{(\alpha)} \neq \emptyset $ then by Lemma 2.6, $ \left| {E_D^{(\alpha )} } \right| = \aleph $ and thus $\left| E \right|= \aleph$, a contradiction.
Thus $E_D^{(\alpha)}=\emptyset$ and $S$ is a $C_\alpha $-sequence, as required.
\end{proof}
Now we consider the possibility of eliminating a non$C_\alpha$-point of a sequence of meromorphic functions. First we have the following definition.
\begin {definition} (cf. \cite [Definition 5.6]{Chuang}, \cite [Definition 3.10]{Applications}, \cite[Definition 8.15]{Chuang}) 
Let $S=\left\{{f_n}\right\}$ be a sequence of meromorphic functions in a domain $D$ and $z_0$ a non$C_0$-point of $S$ in $D$.
Then two cases are possible:\newline
(1) We can extract from the sequence $S$ a subsequence $S'$ of which $z_0$ is a $C_0$-point.\newline
(2) $z_0$ is a non$C_0$-point of every subsequence of the sequence $S$. In this case, we say that $S$ is irreducible with respect to the point $z_0$.
\end{definition}
A sequence is said to be an irreducible sequence in $D$ if $S$ has non$C_0$-points in $D$ with respect to which $S$ is irreducible.
If, in addition, $S$ is irreducible with respect to each of its non$C_0$-point, then $S$ is called completely irreducible sequence.
We have the following theorem.
\begin{theorem} (cf. \cite[Lemmas 5.5, 8.8]{Chuang}) 
Let $S_0  = \left\{ {f_n } \right\}_{n = 1}^\infty  $ be a sequence of meromorphic functions in $D$, and assume that the set $E$ of non$C_0$-points of $S_0$ is at most countable.
Then we can extract from $S_0$ a subsequence $S$ which is either a $C_0$-sequence in $D$ or a completley irreducible sequence in $D$.
\end{theorem}
\begin{proof} (cf. \cite[Proofs of Lemmas 5.5, 8.8]{Chuang})
We assume that $E$ is not finite. The proof for the finite case will be clear after we prove the infinite case. Let $E = \left\{ {z_n } \right\}_{n = 1}^\infty $ be an enumeration of $E$. If $S_0$ is reducible with respect to $z_1$, then there is a subsequence, say $S_1$, such that $z_1$ is a $C_0$-point of $S_1$. If this is not the case we denote $S_1:=S_0$. Suppose we have defined $S_1,...,S_k$.
Then in the $(k+1)$'th step the subsequence $S_{k+1}$ of $S_k$ is determined as follows. If $S_k$ is reducible with respect to $z_{k+1}$ then $S_{k+1}$ is some subsequence of $S_k$, such that $z_{k+1}$ is a $C_0$ sequence of $S_{k+1}$. If this is not the case then $S_{k+1}:=S_k$. Performing this process successively for $k=1,2,...$, we get subsequences $S_1,S_2,...$ such that $S_{k+1}$ is a subsequence of $S_k$. We then define $S$ to be the diagonal sequence corresponding to these sequences, i.e, the first element is the first element of $S_1$, the second element is the second element of $S_2$ and so on. It is clear that $S$ has the required property.
\end{proof}
Observe that if $S_k\neq S_{k+1}$ for every $k \geq 0$, then S is a $C_0$ -sequence in $D$. This process is not unique and we can get various subsequences for $S$, each of them with its set of non$C_0$-points with respect to which it is irreducible. For a sharper result in this issue see \cite{NPZ}.

We can extend  Definition 2.23 as follows.
\begin{definition} (cf.\cite[Definition 8.15]{Chuang}) 
Let $S$ be a sequence of meromorphic functions in $D$ and let $z_0\in D$ be a non$C_\alpha$-point of $S$ for some ordinal number $\alpha$. Then $S$ is said to be \textit{reducible with respect to $z_0$} if there is a subsequence $\hat{S}$ of $S$ such that $z_0$ is a \textit{$C_\alpha$}-point of $\hat{S}$. Otherwise we say that $\hat{S}$ is irreducible with respect to $z_0$.
\end{definition}
In view of Lemmas 2.21, 2.22 and Theorem 2.24, we can state
\begin{corollary} 
Let $S$ be a $C_{\alpha}$-sequence of meromorphic functions in $D$ which is not a $C_\beta$-sequence for any $\beta<\alpha$. Then $S$ has a subsequence $\hat S$ which is either a $C_{\alpha}$-sequence which is not a $C_{\beta}$-sequence for any $\beta<\alpha$ and $\hat S$ is irreducible with respect to each of its non$C_\beta$-points for any $\beta\leq\alpha$, or that $\hat{S}$ is a $C_{\beta _0 } $-sequence for some $\beta_0<\alpha$ and $\hat{S}$ is irreducible with respect to each of its non$C_\beta$-points for any $\beta\leq\beta_0$.
\end{corollary}
We introduce now a generalization of Marty's theorem for $Q_\alpha$-normal families.
We begin with the following definition
\begin{definition} \cite[Definition 8.9]{Chuang} 
Let $F$ be a family of meromorphic functions in a domain $D$ and $z_j$ ($j=1,2,...,n$) a system of points of $D$. We say that the family $F$ satisfies the condition (M) with respect to the system $z_j$ ($j=1,2,\dots,n$), if there exist disks ${\Delta_j}=\overline{\Delta}(z_j,r_j)$ ($j=1,2,\dots,n$) belonging to $D$ and a number $A>0$ such that for each function $f(z)\in F$, we have
$\mathop {\min }\limits_{1 \le j \le n}  {\mathop {\max }\limits_{z \in \Delta_j } \partial \left( {z,f} \right)}  \le A$.
\end{definition}
The following theorem was stated (with somewhat different notation) and proved in \cite[Theorem 8.16]{Chuang} for  $\alpha=m\in \mathbb{N}$. Here we generalize it for every countable ordinal number.
\begin{theorem}
Let $F$ be a family of functions meromorphic in $D$, and let $\alpha$ be a countable ordinal number. Then $F$ is $Q_\alpha$-normal family in $D$ if and only if every set $E \subset D$ such that $E_D^{(\alpha)}=\emptyset$ contains a finite subset $\hat{E}$ such that $F$ satisfies the condition (M) with respect to~$\hat{E}$.
\end{theorem}
We need some preparations before we go to the proof of Theorem~2.28.
\begin{lemma}
Suppose that $E\subset D$ and $z_0\in E_D^{(\alpha)}$ where $\alpha$ is a countable ordinal number. Then there is $ \hat{E}\subset E $ such that $\hat{E}_D^{(\alpha)}=\left\{{z_0}\right\}$.
\end{lemma}
\begin{proof}
We apply transfinite induction. For $\alpha=0$ the lemma is obvious. Assume that the lemma holds for every ordinal number $\beta\in \Gamma(\alpha)$ and we show that it holds also for $\alpha$.
We separate into two cases:\newline
\textbf{Case (A): $\alpha$ has an i.p.} There is a sequence of different points in
$E_D^{(\alpha  - 1)}, \left\{ {z_n } \right\}_{n = 1}^\infty$ such that $z_n\rightarrow z_0$ and also $z_n\neq z_0$ for every $n\geq 1$. We can find positive numbers $\left\{ {r_n } \right\}_{n = 1}^\infty$ such that the collection
$\left\{ {\Delta (z_n ,r_n )} \right\}_{n = 1}^\infty $ composed of strongly disjoint disks, each of which is contained with its closure in $D$. By the induction assumption and Lemma 2.9, there exists for each $n \geq 1$ a set $E_n  \subset \Delta (z_n ,r_n )\cap E$, such that $(E_n)_D^{(\alpha-1)}=\left\{z_n\right\}$. We then set
$\hat{E}=\mathop  \bigcup \limits_{n = 1}^\infty  E_n $
and by Lemma 2.12
$\hat E_D^{(\alpha  - 1)}  = \left\{ z_n \right\}_{n = 1}^\infty   \cup \left\{ {z_0 } \right\}$ and thus $\hat{E}_D^{(\alpha)}=\left\{{z_0}\right\}.$\newline
\textbf{Case (B): $\alpha$ is a limit ordinal number.}
We set an enumeration of $\Gamma{(\alpha)}$, $\Gamma (\alpha ) = \left\{ {\beta _n } \right\}_{n = 1}^\infty$.
We can also choose an increasing subsequence ${\left\{\beta _{n_k }\right\}}_{k=1}^\infty$, $\beta _{n_k }  \nearrow \alpha $ (i.e., for every $\beta < \alpha$ there is $k_0=k_0(\beta)$ such that $\beta_{n_k} >\beta$ for $k > k_0$). By Definition 2.1, there are points $\left\{{z_n}\right\}$ in $D$ such that $ \left|z_n-z_0\right|\rightarrow 0 $ and $z_n \in E_D^{(\beta_n)}$. As in Case (A), we can find $\left\{ {r_n } \right\}_{n = 1}^\infty $ such that $\left\{ {\Delta (z_n ,r_n )} \right\}_{n = 1}^\infty$ are strongly disjoint in $D$. By the induction assumption, for every $n\geq 1$ there is $E_n\subset E\cap \Delta(z_n,r_n)$ such that
$(E_n )_D^{\left( {\beta _n } \right)}  = \left\{ {z_n } \right\}_{n = 1}^\infty  $.
Again, the set
$\hat{E}=\mathop  \cup \limits_{n = 1}^\infty  E_n$
satisfies $\hat{E}_D^{(\alpha)}=\left\{{z_0}\right\}$, as desired.
\end{proof}
We also need the following definition.
\begin{definition} \cite[Definition 3.4]{Transfinite} 
Let $S=\left\{{f_n}\right\}$ be a sequence of meromorphic functions on a domain $D$ and let $E\subset D$. Then $E$ and $S$ are said to satisfy Zalcman's condition with respect to $D$ if for every $\zeta\in E $ there exist $z_n,_\zeta\rightarrow \zeta$, $z_n,_\zeta \in D$; $\rho_{n,\zeta}\rightarrow 0^{+}$, such that
$f_n \left( {z_{n,\zeta }  + \rho _{n,\zeta } \zeta } \right)\mathop  \Rightarrow \limits^\chi  g_\zeta  \left( z \right)$ on $\mathbb{C}$, where $g_\zeta$ is a nonconstant meromorphic function on $\mathbb{C}$.
Observe that in the setting of Definition 2.30, every point of $E$ is a non$C_0$-point of $S$, with respect to which $S$ is irreducible.

\end{definition}

\begin{theorem} \cite[Theorem 3.5]{Transfinite} 
Let $\alpha$ be a countable ordinal number. Then a family $F$ of meromorphic functions in the domain $D$ is not $Q_\alpha$-normal if and only if there is $E\subset D$ satisfying $E_D^{(\alpha)}\neq \emptyset $, and a sequence $S$ of functions from $F$ such that $S$ and $E$ satisfy Zalcman's condition with respect to $D$.
\end{theorem}
\begin{proof}[Proof of Theorem 2.28]
Assume first that $F$ is a $Q_\alpha$-normal family, and let $E\subset D$ be such that $E_D^{(\alpha)}\neq \emptyset$. By Lemma 2.29 there is $\hat{E}\subset E$ such that $\left|\hat{E}_D^{(\alpha)}\right|=1$. We can also assume, without loss of generality, that $dist(\hat{E},\partial{D}) > 1$.
Let the enumeraion of $\hat E$ be $\hat E = \left\{ {a_n } \right\}_{n = 1}^\infty$. If on the contrary, the family $F$ does not satisfy the condition (M) with respcet to any finite subset of $\hat E$, then there is a sequence in $F$, $S = \left\{ {f_n } \right\}_{n = 1}^\infty$ such that for every $n \geq 1$ we have
 \begin{equation}
\mathop {\max f_i^\#  (z)}\limits_{z \in \bar \Delta \left( {a_i ,\frac{1}{n}} \right)}  > n
\label{eq:2.5}
\end{equation}
for $1\leq i \leq n $.
By Marty's Theorem we get by (\ref{eq:2.5}) that every point in $\hat{E}$ is a non$C_0$-point of $S$, with respect to which $S$ is irreducible. Since $\hat{E}_D^{(\alpha)}\neq \emptyset$, we have a contradiction to the $Q_\alpha$-normality of $F$.

In the other direction, suppose that the family $F$ is not a $Q_\alpha$-normal family in $D$. Then by Theorem 2.31  there is a sequence $S=\left\{{f_n}\right\}_{n=1}^\infty$ in $F$, and a set $E\subset D$, $E_D^{(\alpha)}\neq \emptyset$ such that $S$ and $E$ satisfy Zalcman's condition with respect to $D$. Without loss of generality, we can assume that $d=dist(E,\partial D)>0$.
Let $\tilde{E}=\left\{{z_1,...,z_m}\right\}$ be an arbitrary finite subset of $E$. Then for every $r<d$ we have for every $1 \leq i \leq m$,
$\mathop {\max f_n^\#  (z)}\limits_{z \in \bar \Delta \left( {z_i ,r} \right)}  \to \infty $.
Thus $F$ does not satisfy the condition (M) with respect to $\hat E$, a contradiction. We deduce that $F$ is a $Q_\alpha$-normal family in D, as desired.
\end{proof}
\begin{remark} 
Denote by $\alpha_{\aleph}$ the first noncountable ordinal number. Lemma 2.22 makes it superfluous to consider $Q_\alpha$-normal  for $\alpha > \alpha_{\aleph}$. Indeed, by applying transfinite induction and using a similar argument to that in the proof of Lemma 2.21, we can deduce that if $S$ is a $C_\beta$-sequence for some ordinal number $\beta$, then the set $E$ of its non$C_0$-points of $S$ is at most enumerable, and then by the other direction of Lemma~2.22, $S$ is already a $C_{\alpha}$-sequence for some countable ordinal number $\alpha$. The value of $\alpha$ can be in general large as we like, and thus by Definition~2.17 we can say that any $Q_\alpha$-normal family is already a $Q_{\alpha_{\aleph}}$-normal family, but in general not a $Q_\alpha$-normal family for any $\alpha<\alpha_{\aleph}$.
\end{remark}
It is necessary to mention that for every countable ordinal number $\alpha$, there is a $C_\alpha$-sequence in $D$ which is not a $C_\beta$-sequence for every $\beta<\alpha$ \cite[Theorem 3.1]{Transfinite}. The way to construct such a family is to find a set $E \subset D$ such that $E_D^{(\beta)}\neq \emptyset$ for every $\beta<\alpha$ and $E_D^{(\alpha)}=\emptyset$. (This is done by applying transfinite induction.) Then one only has to use the following lemma and Corolloary 2.13.
\begin{lemma} (\cite[Lemma 8.4]{Chuang}, \cite[Lemma 3.2]{Transfinite}) 
Let $E\subset$ $D$ satisfy $E\cap E_D^{1}=\emptyset$. Then there is a sequence $S$ of polynomials, such that $E\cup E_D^{1}$ is exactly the set of non$C_0$-points of $S$, with respect to each of which $S$ is irreducible.
\end{lemma}
We will also use this principle in the proof of Theorem 3.1 (Main Theorem); See Lemma 2.36.
Now we define the notion of order of $Q_\alpha$-normal family.
\begin{definition} (cf. \cite[Definition 8.7]{Chuang}, \cite[Definition 3.29]{Applications})
Let $\alpha$ be a countable ordinal number with i.p. Let $F$ be a family of meromorphic functions in a domain $D$ and $\nu\geq 1$ an integer. We say that $F$ is $Q_\alpha$ -normal of order at most $\nu$ in $D$ if from every sequence of functions of the family $F$ we can extract a subsequence which is a $C_\alpha$-sequence in $D$ and has at most $\nu$ non$C_{\alpha-1}$-points in $D$. In particular, when $\nu=0$, F is  $Q_{\alpha-1}$-normal in $D$.
If $F$ is a $Q_\alpha$-normal family of order at most $\nu\geq 1$ in $D$ but not a $Q_m$ normal family of order at most $\nu-1$ in $D$, we say that $F$ is a $Q_\alpha$-normal family of exact order $\nu$ on $D$. If $F$ is a $Q_\alpha$-normal family but not a $Q_\alpha$-normal family of order at most $\nu$ for any $\nu\geq 1$, we say that $F$ is a is a $Q_\alpha$-normal family of infinite order in $D$.
\end{definition}
\begin{remark} 
There is probably no reason to define a $Q_\alpha$-normal family of order at most $\nu$  (where $\nu<\infty$) when $\alpha$ is a limit ordinal number. Indeed a reasonable definition, consistent with Definition 2.34, would be that a family $F$ is $Q_\alpha$-normal of order at most $\nu$, if every sequence  $S$ in $F$ has a subsequence $\hat{S}$ with at most $\nu$ points in $D$, which are non$C_\beta$-points for every $\beta<\alpha$. But by Definition 2.15 (for limit ordinal numbers) there are no such points in a $C_\alpha$-sequence ($\hat{S}$ is of course a $C_\alpha$-sequence by that condition).
Another possibility might be that for every sequence $S$ in $F$, there is some $\beta<\alpha$ and a subsequence $\hat{S}$ with at most $\nu$ non$C_\beta$-points. But then such a finite number of points are all $C_{{\beta}'}$-points for some $\beta<{\beta}'<\alpha$, and then by the same definition $F$ is a $Q_\alpha$-normal family of order $0$.
To summarize, the problem in defining this notion for limit ordinal numbers for every sequence in such family should have a subsequence with at most $\nu$ `problematic' points. But all these problematic points are $C_{\beta _0 }$-points for some $\beta_0<\alpha$,and being $\beta_0$ `very far' from $\alpha$ (since $\alpha$ is a limit ordinal number) makes difficult to define this `problematic nature' of arbitrary points in a satisfactry way.
In the case $\nu=\infty$, i.e., defining $Q_\alpha$-normal family of infinite order, the second suggestion for such a definition make sense. Indeed, if there is a subsequence $\hat{S}$ with infinitely many points $\left\{ {z_n } \right\}_{n = 1}^\infty   \subset D $, all of which are non$C_\beta $-points for some $\beta<\alpha$, then we can assume that $\hat{S}$ is irreducible with respect to each $z_n$ (otherwisw we apply the diagonal process as in Theorem 2.24). If there is always some ${\beta}'<\alpha$ such that each $z_n$ is a $C_{{\beta}'}$-point, then again there is no purpose to the definition.
We will get a `real' $Q_{\alpha}$-normal family of order $\infty$ if there is a case where such ${\beta}'$ does not exist. In such a case let $z_0$ be an accumulation points of $\left\{ {z_n } \right\}_{n = 1}^\infty$. It is impossible that $z_0  \in D$, because then $z_0$ would be a non$C_\alpha$-point that $\hat S$ is irreducible with respect to it. Thus $z_0=\infty$ or $z_0 \in \partial D$. However, our main theorem, Theorem 1, deals with the case  where $D=\Delta$, and then the possibility $z_0=\infty$ is excluded and the possibility $\left| z_0 \right|=1$ would lead to a contradiction, because of the special nature of the famiy $F(f)$ (See section 4.1). Thus we have not defined a $Q_{\alpha}$-normal family of infinite order in the case where $\alpha$ is a limit ordinal number.
We add that, in this case, it is possible by Lemma 2.33 to create families that are $Q_\alpha$-normal of order at most $\nu$ (finite or infinite) according to the second suggested definition.
\end{remark}
The following Lemma is essential for the proof of our main theorem.
\begin{lemma} (cf. \cite [Lemma 3.3]{Transfinite}, \cite [Lemma 3.3]{Generating}) 
Let $\alpha$ be a finite or countable ordinal number and let $\Gamma$ be any compact arc of a circle. Then if $\alpha-1$ exists and $ 1\leq\nu < \infty$ is an integer, there exists a set $ E=E(\alpha,\nu,\Gamma)\subset \Gamma$ such that $E\cap E_\Gamma^{1}=\emptyset$ and $\left| {E_\Gamma ^{(\alpha  - 1)} } \right| = \nu $ (and thus $E_\Gamma^{(\alpha)}=\emptyset$).
If $\alpha$ is a limit ordinal number, then there is a set $ E=E(\alpha,1,\Gamma)\subset \Gamma $ satisfying $\left| {E_\Gamma ^{(\alpha )} } \right| = 1$ and $E\cap E_\Gamma^{(1)}=\emptyset$.
\end{lemma}
\begin{proof}
Without loss of generality, we can assume that $\Gamma$ is on the unit circle. We will apply transfinite induction and assume first that $\nu=1$ (for the case that $\alpha-1$ exists). In case $\alpha=1$ we take $E(1,1,\Gamma)=\left\{{z_0}\right\}$ where $z_0 \in \Gamma$. In case $\alpha=2$ let $z_0=e^{i\theta_0} \in \Gamma $ and we take a sequence $\left\{ {\theta _n } \right\}_{n = 1}^\infty $ such that (without loss of generality) $\theta_n \nearrow \theta_0$ and (2.1)
\begin {equation}
 z_n  = e^{i\theta _n }  \in \Gamma   \quad for\quad every  \quad n\geq 1.
\end {equation}
Then set $E(2,1,\Gamma ) = \left\{ {z_n } \right\}_{n = 1}^\infty$.
Now let $\alpha$ be any countable (or finite) ordinal number. We separate into two cases:\newline
\textbf{Case (A) $\alpha-1$ exists.} Enclose each $z_n$ from (2.6) in a disc $\Delta_n=\Delta(z_n,r_n)$ such that $\left\{ {\Delta _n } \right\}_{n = 1}^\infty$ are strongly disjoint. We separate now into two subcases.\newline
\textbf{Case $(A_1)$ $\alpha-2$ exists}. For every $n\geq 1$ denote
\begin{equation}
\Gamma_n=\Gamma\cap {\overline{\Delta}_n}.
\end{equation}
By the induction assumption there exists, for every $n \geq 1$, a set $E_n=E(\alpha-1,1,\Gamma_n)\subset\Gamma_n$, such that
 $\left| {(E_n)}_\Gamma^{(\alpha  - 2)} \right| = 1$ and $E_n  \cap {(E_n)}_\Gamma ^{(1)}  = \emptyset$. Set $E = E(\alpha ,1) = \bigcup\limits_{n = 1}^\infty  {E_n } $. Then by Lemma 2.12 $E$ has the required property.\newline
\textbf{Case $(A_2)$ $\alpha-1$ is a limit ordinal number.} We first arrange $\Gamma(\alpha-1)$ in an enumeration, $\Gamma (\alpha  - 1) = \left\{ {\alpha _n } \right\}_{n = 1}^\infty  $.
For every $\nu \geq 1$, there is, by the induction assumption, a set $E_n=E(\alpha_n,1,\Gamma_n) \subset \Gamma_n$, $\Gamma_n$ as in (2.7). Again we define $E = \bigcup\limits_{n = 1}^\infty  {E_n } $ and it has the required property.\newline
\textbf{Case (B): $\alpha$ is a limit ordinal number.} We write $\Gamma (\alpha ) = \left\{ {\alpha _n } \right\}_{n = 1}^\infty$ and by the induction assumption there exists for every $n \geq 1$ a set $E_n (\alpha _n ,1,\Gamma _n ) \subset \Gamma_n $. Again set $E = \mathop  \bigcup \limits_{n = 1}^\infty  E_n $ and it works.
In the case where $1 \leq \nu < \infty$ and $\alpha-1$ exists, we take $\nu$ strongly disjoint compact subarcs, $\Gamma_1,\Gamma_2,...,\Gamma_\nu$, and for each $1 \le j \le \nu $ there is a set $E_j=E(\alpha,1,\Gamma_j) \subset \Gamma_j$ and we take $E = \mathop  \bigcup \limits_{j = 1}^\nu  E_j $
\end{proof}
We now state our main theorem.
\begin{thm1*}
Let $\alpha$ be a countable (or finite) ordinal number. Then there is an entire function $f(z)=f_\alpha(z)$ such that the family $F(f)$ is a $Q_\alpha$-normal family on $\Delta$, but not a $Q_\beta$-normal family for any $\beta<\alpha$. Morever, if $\alpha-1$ exists and $\nu\geq 1$ is an integer or $\nu = \infty$, then there is an entire function $f(z)=f_{\alpha,\nu}(z)$ such that $F(f)$ is $Q_\alpha$-normal of exact order $\nu$ (or of infinite order if $\nu = \infty$)
\end{thm1*}

The proof of Theorem 1 is devided into 3 parts. The first part is proved for the case where $\alpha-1$ exists.
\section{Proof of Theorem 1 for $\alpha$ having i.p and $\nu<\infty$}
The proof of this case of the theorem, as the proof of the other cases, goes by constructing the required entire function $f=f_{\alpha,\nu}$.
We assume that $\alpha-1$ exists and that $1 \leq \nu< \infty$ is an integer.
\begin{proof}
Let $E=E(\alpha,\nu)\subset\Gamma:=\partial\Delta$ be as guaranteed by Lemma 2.36. $E$ is enumerable by Lemma 2.7 and satisfies
\begin{equation}
 E\cap E_\Gamma^{1}=\emptyset , \left|E_\Gamma^{(\alpha-1)}\right|=\nu.
 \end{equation}
 Let $E=\left\{{c_n}\right\}_{n=1}^\infty$ be an enumeration of E. Define now an increasing sequence of positive numbers that satisfies the following 3 conditions:
\begin{equation}
\left( {\frac{{a_n }}{{a_{n - 1} }}} \right)^{\frac{1}{n}} \mathop  \to \limits_{n \to \infty } \infty \quad.
\end{equation}
For large enough $n$
\begin{equation}
a_n  \ge \frac{1}{{1 - \left( {1 - \frac{1}{{2^{n + 1} }}} \right)^{\frac{1}{{n + 1}}} }} \quad.
\end{equation}
For every $n\geq1$
\begin{equation}
a_{n + 2}  \ge a_{n + 1} a_n \quad.
\end{equation}
It is not hard to see that in fact condition (3.4) implies conditions (3.2) and (3.3). One way to get such a sequence $\left\{{a_n}\right\}$ is by letting $a_n  = e^{\alpha _n } $ where $\left\{{\alpha_n}\right\}_{n=1}^\infty$ is a Fibonacci sequence of positive numbers. Now we define a sequence $\left\{{b_l}\right\}_{l=1}^\infty$, to be the (simple) zeros of $f(z)$. For convenience and explantory purposes, we introduce $\left\{{b_l}\right\}_{l=1}^\infty$ in rows 
from the top down, $R_1,R_2,...,R_k,..$ where the order in each row is from left to right, as follows:
\begin{equation}
\begin{array}{l}
 R_1 :a_1 c_1 ; \\
 R_2 :a_2 c_1, a_2 c_2 ; \\
 \,\,\,\,\,\,\,\,\,\,\, \vdots  \\
 R_n ;a_n c_1, a_n c_2, ...,a_n c_n ; \\
  \,\,\,\,\,\,\,\,\,\,\, \vdots  \\
 \end{array}
\end{equation}
 For example $b_1=a_1c_1$, $b_3=a_2c_2$, $b_7=a_4c_1$, and so on.
 For each $l\geq 1$ define $S(l)$ to be the natural number $n$, where $s(l)=n$ indicates that $b_l\in R_n$ (or that $\left|b_l\right|=a_n$). For example, $S(7)=S(8)=S(9)=S(10)=4$.
 Then
\begin{equation}
\left\{ {l:s(l) = n} \right\} = n.
\end{equation}
It is easy to see (by (3.2) for example) that
$\sum\limits_1^\infty  {\frac{n}{{a_n }}}  < \infty$; and thus by (3.6), the factorization
\begin{equation}
f(z): = \prod\limits_{l = 1}^\infty  {\left( {1 - \frac{z}{{b_l }}} \right)}
\end{equation}
defines an entire funcion, vanishing exactly at $z=b_l$, $l\geq 1$.

Now we are going to prove that $f(z)$ has the required property. In what follows we will consider only sequences
$\left\{ {f_{j_k } } \right\}_{k = 1}^\infty   \subset F(f)$ where $j_k \to \infty$; otherwise we can pick a normal subsequence. First we would like to limit the direction (or ray from the origin) on which non$C_0$-points of some sequence $\left\{ {f_{j_k } } \right\}$ can be.
\begin{lemma}
Let $\eta _0  \in \Delta \backslash \{ L(\arg z):z \in \bar E\}$.
Then $\eta_0$ is a $C_0$-point of $\left\{ {f_n (z)} \right\}_{n = 1}^\infty$  and  $f_n(z)\Rightarrow\infty$ on some neighboerhood of $\eta_0$.
\end{lemma} 
\begin{proof}
By assumption there is some $\alpha_0>0$ such that
\begin{equation}
S(\arg \eta _0 ,\alpha _0 ) \cap \{ b_l :l \ge 1\}  = \emptyset
\end{equation}
We shall show that $f(z)\mathop  \Rightarrow \limits_{z \to \infty } \infty $
in $S(\arg \eta _0 ,\frac{{\alpha _0 }}{2})$, and since $f_n(z)$ attains in a small enough disk $\Delta(\eta_0,r)$ the values that $f(z)$ attains in $\Delta (n\eta _0 ,nr) \subset S(\arg \eta _0 ,\frac{{\alpha _0 }}{2})$, this will prove the lemma. So let
$z \in S(\arg \eta _0 ,\frac{{\alpha _0 }}{2})$
be with large enough modul such that there exists a (unique) $n$ such that
\begin{equation}
a_1<\left|z \right|\leq a_{n+1}.
\end{equation}
We estimate the factors of the product (3.7). All estimations are valid for large enough $n$ (or large enough $l$)\newline
\textbf{(A) $l:S(l)\leq {n-1}$}.

By (3.2) and (3.9)
\begin{align*}
\left| {1 - \frac{z}{{b_l }}} \right| \ge \left| {\frac{z}{{b_l }}} \right| - 1 \ge \frac{{a_n }}{{a_{n - 1} }} - 1 > 2.
\end{align*}
Thus
\begin{align*}
\big| {\prod\limits_{l:s(l) \le n - 1} {\left( {1 - \frac{z}{{b_l }}} \right)} } \big| \ge 2^{1 + 2 + ... + n - 1}  = 2^{\frac{{n(n - 1)}}{2}}.
\end{align*}
\textbf{(B) $l: n\leq S(l)\leq n+2$}.\newline
In this case, by (3.8)
$\big| {\frac{z}{{b_l }} - 1} \big| \ge \sin \frac{{\alpha _0 }}{2} > 0$
and then
\begin{align*}
\big| {\prod\limits_{l:s(l) = n,n + 1,n + 2}^{} {\left( {1 - \frac{z}{{b_l }}} \right)} } \big| \ge \left( {\sin \frac{{\alpha _0 }}{2}} \right)^{n + (n + 1) + n + 2}  = \left( {\sin \frac{{\alpha _0 }}{2}} \right)^{3(n + 1)}.
\end{align*}
\textbf{(C) $l:S(l) \ge n + 3$}.\newline
We can write for such $l$, $s(l)=n+k$, $k\geq 3$, and then, by the definition of $S(l)$, (3.4) and (then) (3.3), we have
\begin{align*}
\big| {\prod\limits_{l:s(l) \ge n + 3}^{} {\left( {1 - \frac{z}{{b_l }}} \right)} } \big| \ge \prod\limits_{j = n + 3}^\infty  {\big| {\prod\limits_{l:s(l) = j} {\left( {1 - \frac{z}{{b_l }}} \right)} } \big|}  \ge \prod\limits_{j = n + 3}^\infty  {\left[ {\left( {1 - \frac{1}{{2^j }}} \right)^{\frac{1}{j}} } \right]^j }
\end{align*}
\begin{align*}
= \prod\limits_{j = n + 3}^\infty  {\left( {1 - \frac{1}{{2^j }}} \right) \ge k_0 : = } \prod\limits_{j = 1}^\infty  {\left( {1 - \frac{1}{{2^j }}} \right) > 0} .
\end{align*}

 By collecting the results of (A), (B) and (C), we get that
$
\left| {f(z)} \right| \ge 2^{\frac{{n(n - 1)}}{2}} \left( {\sin \frac{{\alpha _0 }}{2}} \right)^{3(n + 1)}  \cdot k_0  \ge \left( {\frac{3}{2}} \right)^{\frac{{n(n - 1)}}{2}}  \to \infty
$
and the lemma is proved.
\end{proof}

The next step is to show that when we have a non$C_0$-point of some $\left\{{{f_j}_k}\right\}\subset F(f)$, then we can pick some `controllable' subsequence of $\left\{{{f_j}_k}\right\}$.
\begin{lemma}
Let $\left\{ {f_{j_k } } \right\}_{k = 1}^\infty  $ be a sequence in $F(f)$, and suppose that an element $0\neq\eta_0\in \Delta$ is a non$C_0$-point of $\left\{{f_{j_k}}\right\}$. For large enough $k$ we can write $a_{n_k }  < j_k \left| {\eta _0 } \right| \le a_{n_k  + 1} $, where $\eta _k  = n(\left| {\eta _0 } \right| \cdot j_k )$.
Then there is a subequence $\left\{{j_{k_p}}\right\}_{p=1}^\infty$ of $\left\{{j_k}\right\}$ that fulfills one of the following options: \newline
\begin{equation}
(1)\quad
\frac{{j_{k_p } \left| {\eta _0 } \right|}}{{a_{n_{k_p } } }}\mathop  \to \limits_{p \to \infty } 1^ +
\quad or \quad (2) \quad
\frac{{j_{k_p } \left| {\eta _0 } \right|}}{{a_{n_{k_p }  + 1} }}\mathop  \to \limits_{p \to \infty } 1^ -.
\end{equation}
In addition, $\left\{{f_{j_{k_{p}}}}\right\}$ has no non$C_0$-points outside $\left\{{0}\right\}\cup\left\{{z:\left|z\right|=\eta_0}\right\}$.
\end{lemma}
\begin {proof}
If none of the two possibilities in (3.10) holds, then there exists $0<\delta<1$, such that for large enough $k$
\begin{equation}
\frac{{j_{k_{} } \left| {\eta _0 } \right|}}{{\left| {a_{n_k  + 1} } \right|}} < 1 - \delta
\quad and \quad
  \frac{{j_{k_{} } \left| {\eta _0 } \right|}}{{a_{n_k } }} > 1 + \delta .
\end{equation}
Suppose, without loss of generality, that (3.11) holds for every $k\geq 1$ and that $\delta$ is so small such that
$\Delta  \supset \bar \Delta (\eta _0 ,\left| {\frac{{\eta _0 }}{2}} \right|\delta ).$
Now if
$z \in \Delta (\eta _0 ,\frac{{\left| {\eta _0 } \right|\delta }}{2})$, then
\begin{equation}
j_k \left| {\eta _0 } \right|\Big(1 - \frac{\delta }{2}\Big) < j_k \left| z \right| < j_k \left| {\eta _0 } \right|\Big(1 + \frac{\delta }{2}\Big).
\end{equation}
For such $z$ we estimate the product
$f_{j_k } (z) = \prod\limits_{l = 1}^\infty  {(1 - \frac{{j_k z}}{{b_l }})}$,
doing it for large enough $k$. As in the proof of Lemma 3.1, also here we separate the range of values of $l$ ($l\geq1$).\newline
\textbf{(A) $l:s(l)\leq\ n_k-1$}.\newline
By (3.11) and (3.12) \newline
$\left| {\frac{{j_k z}}{{b_l }}} \right| \ge \frac{{j_k \left| {\eta _0 } \right|(1 - \frac{\delta }{2})}}{{a_{n _k  - 1} }} >   (1 + \delta )(1 - \frac{\delta }{2})\frac{{a_{n _k } }}{{a_{n_k  - 1} }}$. Evidently for such $l$, $
\left| {1 - \frac{{j_k z}}{{b_l }}} \right| \ge 2$, and
$
\Big| {\prod\limits_{l:s(l) \le n_k - 1}^{} ( 1 - \frac{{j_k z}}{{b_l }})} \Big| \ge 2^{\frac{{(n_k  - 1)n_k }}{2}}.
$\newline
\textbf{(B)} $l: n_k \leq s(l)\leq n_k+2$.\newline
 First we observe that if $s(l)=n_k$, then
$
\left| {\frac{{j_k z}}{{b_l }}} \right| > \frac{{j_k \left| {\eta _0 } \right|(1 - \frac{\delta }{2})}}{{\left| {b_{l_k } } \right|}} > (1 + \delta )(1 - \frac{\delta }{2}) > 1
$ and
\begin{equation}
\left| {1 - \frac{{j_k z}}{{b_l }}} \right| > (1 + \delta )\Big(1 - \frac{\delta }{2}\Big) - 1 = \frac{{\delta (1 - \delta )}}{2}.
\end{equation}
On the other hand, if $s(l)=n_k+1$ or if $s(l)=n_k+2$, then
$
\left| {\frac{{j_k z}}{{b_l }}} \right| < \frac{{j_k \left| {\eta _0 } \right|(1 + \frac{\delta }{2})}}{{\left| {b_{l_{} } } \right|}} \le \frac{{j_k \left| {\eta _0 } \right|(1 + \frac{\delta }{2})}}{{a_{n_k  + 1} }} < (1 - \delta )(1 + \frac{\delta }{2}) < 1
$
and thus
\begin{equation}
\left| {1 - \frac{{j_k z}}{{b_l }}} \right| > 1 - (1 - \frac{\delta }{2})(1 + \frac{\delta }{2}) = \frac{{\delta (1 + \delta )}}{2}.
\end{equation}
By (3.13) and (3.14) we have
$
\Big| {\prod\limits_{l:s(l) = n _k ,n_k  + 1,n_k  + 2}^{} ( 1 - \frac{{j_k z}}{{b_l }})} \Big| \ge \left[ {\frac{{\delta (1 - \delta )}}{2}} \right]^{3(n _k  + 1)}.
$\newline
\textbf{(C) $l:s(l)\geq n_k+3$}.\newline
 We can write $s(l)=n_k+p$ with $p \geq 3$. Then by (3.11) and (3.12) and (3.3), we get
\begin{align*}
\left| {1 - \frac{{j_k z}}{{b_l }}} \right| > 1 - \frac{{\left| {j_k z} \right|}}{{a_{n_k  + p} }} > 1 - \frac{{(1 - \delta )(1 + \delta /2)}}{{a_{n_k  + p - 1} }}
\end{align*}
\begin{align*}
\ge 1 - \frac{1}{{a_{n_k  + p - 2} }} \ge \left( {1 - \frac{1}{{2^{n_k  + p - 1} }}} \right)^{\frac{1}{{n_k  + p - 1}}}.
\end{align*}
Hence we have

\begin{align*}
\Big| {\prod\limits_{l:s(l) \ge n _k  + 3}^{} ( 1 - \frac{{j_k z}}{{b_l }})} \Big| \ge \Big| {\prod\limits_{p = 3}^\infty  ({\prod\limits_{l:s(l) = n _k  + p} {(1 - \frac{{j_k z}}{{b_l }})} } }) \Big|
\end{align*}
\begin{align*}
 \ge \Big| {\prod\limits_{p = 3}^\infty  {(1 - \frac{1}{{2^{n _k  + p - 1} }})^{\frac{{n _{k }  + p}}{{n _k  + p - 1}}} } } \Big| \ge \prod\limits_{t = 1}^\infty  {(1 - \frac{1}{{2^t }}} )^2  = k_0^2  > 0. \newline
\end{align*}
From the estimates in (A), (B) and (C), we deduce similarly as in the proof of Lemma 3.1 that $f_{j_k } \left( z \right)\Rightarrow \infty $ on
$\Delta (\eta _0 ,\frac{{\left| {\eta _0 } \right|}}{2}\delta )$, and this is a contradiction for $\eta_0$ to be a non$C_0$-point of $\left\{{f_{j_k}}\right\}$.
Now, if $0 < r \ne \left| {\eta _0 } \right|$, then $zr$ is a $C_0$-point of $\left\{f_{j_{k_p}}\right\}$ for every $z\in\partial\Delta$. This follows from (3.10), since
$\frac{{\left| z \right|}}{{\left| {\eta _0 } \right|}} \ne 0,1,\infty $, and $\frac{{a_{n _{k_p }  + 1} }}{{a_{n _{k_p } } }}\mathop  \to \limits_{p \to \infty } \infty$. This completes the proof of Lemma 3.2 .
\end{proof}
Lemma 3.2 implies that $F(f)$ is a $Q_\alpha$-normal family in $\Delta$ (in fact in $\mathbb{C}$). The next lemma shows that $F(f)$ is not a $Q_\beta$-normal family in $\Delta$ for any $\beta < \alpha$. This will complete the proof of Theorem 1 for the case where $\alpha-1$ exists and $\nu < \infty$.
\begin{lemma}
For every $0<r<1$ there is a sequence $\left\{ {f_{j_k } } \right\}_{k = 1}^\infty $, $j_k=j_k(r)$ such that the set of non$C_0$-points of every subsequence of $\left\{ {f_{j_k } } \right\}_{k = 1}^\infty$ is exactly $\left\{{0}\right\}\cup r\overline{E}$ (i.e., $\left\{ {f_{j_k } } \right\}_{k = 1}^\infty$ is irreducible with respect to each point of this set).
\end{lemma}
\begin{proof}
Set
\begin{equation}
j_k : = \left[ {\frac{{a_k }}{r} + 1} \right].
\end{equation}
We have
$a_k  < j_k r < a_{k + 1} $
, and also
$\frac{{j_k r}}{{a_k }} \to 1^+$.
This means that option (1) in (3.10) holds, and thus by Lemma 3.2 every $z'$ such that $\left|z'\right|\neq 0,r$ is a $C_0$-point of $\left\{{f_{j_k}}\right\}$. We claim that for every $z\in E$, $rz$ is a non$C_0$-point of any subsequence of $\left\{{f_{j_k}}\right\}$. By Lemma 3.1 it is enough to show that for every $\delta>0$, $f_{j_k}$ has zero in $\Delta(rz,r\delta)$ for large enough $k$. That is, we have to show  that $\Delta(j_krz,j_kr\delta)$ contains a zero of $f(z)$, for large enough $k$. Indeed, by (3.15)  $\zeta \in \Delta(j_krz,j_kr\delta)$ if and only if
$\Big| {\frac{{\zeta r}}{{a_k }}\cdot\frac{{a_k /r}}{{\left[ {\frac{{a_k }}{r}} \right] + 1}} - rz} \Big| < r\delta $.
We have that $z=c_n$ for some $n\geq1$ and by (3.5) $\zeta_k:=a_kc_n$ is a zero of $f(z)$ for $k\geq n$.

Now $
\Big| {\frac{{j_k r}}{{a_k }} \cdot \frac{{a_k /r}}{{\left[ {\frac{{a_k }}{r} + 1} \right]}} - rz} \Big| = \Big| {c_n r\frac{{a_k r}}{{\left[ {a_k /r + 1} \right]}} - rz} \Big|\mathop  \to \limits_{k \to \infty } 0$
, and so it is less than $r \delta$ for large enough $k$, as required. Hence the set of non$C_0$-points of every subsequence of $\left\{{f_{j_k}}\right\}$ is $A_r:=r\hat{E}\cup{0}$, and by (3.1)
$\left| {\left( {A_r } \right)_\Delta ^{(\alpha  - 1)} } \right| = \nu$, and so we get a sequence $\left\{f_{j_k}\right\}$ which by Lemma 2.20 is by itself a $Q_\alpha$-normal family of exact order $\nu$, and so is $F(f)$.
The proof of Theorem~1 for $\alpha$ with i.p and $\nu<\infty$ is completed.
\end {proof}
We would like to conclude this section with some observations on the nature of non$C_0$-sequences in $F(f)$ for the function $f(z)$ we have constructed.
\begin{lemma}
Let $r>0$ and $z_0\in\bar{E}$. Then $rz_0$ is a non$C_0$-point of a sequence $\left\{{f_{j_k}}\right\}$ in $F(f)$ if and only if there exists a subsequence $\left\{{f_{j_{k_{p}}}}\right\}$ such that for every $\delta>0$ $f_{j_{k_p}}(z)$ has zero in $\Delta(rz_0,\delta)$ for large enough~$p$.
\end{lemma}
\begin{proof}
If $rz_0$ were a $C_0$-point of $\left\{{f_{j_k}}\right\}$, then for small enough $\delta>0$ we would have $f_{j_k}(z)\Rightarrow\infty$ on $\Delta(rz_0,\delta)$ and then for $k \geq k_0$ $f_{j_k}(z)\neq 0$ on $\Delta(rz_0,\delta)$, a contradiction. In the other direction, assume that $rz_0$ is a non$C_0$-point of $\left\{{f_{j_k}(z)}\right\}$.
By Zalcman's Lemma, $rz_0$ is a non$C_0$-point of some sumsequence of $\left\{f_{j_k}\right\}$, say $S_p  = \left\{ {f_{j_{k_p } } } \right\}_{p = 1}^\infty $, such that $S_p$ is irreducible with respect to $rz_0$.
If to the contrary there were a subsequence $S_l=\left\{{f_{j_{k_{p_l}}}}\right\}$ such that $f_{j_{k_{p_l}}}(z)\neq 0$ in $\Delta(rz_0,\delta)$, we would get a contradiction. Indeed, because of the $Q_\alpha$-normality of $F(f)$ there is a subsequence $S_m=\left\{f_{j_{k_{p_{l_m}}}}\right\}$ which is a $C_\alpha$-sequence in $\Delta$ and thus has at most countable set of non$C_0$-points. Thus there is some $0< \delta' <\delta$ such that all the points of the set $\left\{{z:\left|z-z_0\right|=\delta'}\right\}$ are $C_0$-points of $S_m$. Then by the minimum principle we get that $f_{j_{k_{p_{l_m } } } } (z) \Rightarrow \infty $ on $\Delta(rz_0,\delta')$, and this contadicts that $\left\{{f_{j_{k_p}}}\right\}$ is irreducible with respect to $rz_0$.
\end{proof}
\begin{lemma}
If $rz^*$, $z^* \in \overline{E}$ is a non$C_0$-point of $\left\{ {f_{j_k } } \right\}_{k = 1}^\infty $ in $F(f)$, then $rz$ is a non$C_0$-point of $\left\{ {f_{j_k } } \right\}$ for every $z\in \overline{E}$.
\end{lemma}
\begin{proof}
By Lemma 3.2 (with its notations) there is a subsequence $S=\left\{{f_{j_{k_{p}}}}\right\}$ such that
$\frac{{j_{k_p } }r} {{a_{n_{k_p } } }} \to 1^ +  $
(without loss of generality, we assume that possibility 1 in (3.10) holds). Thus for every $z\in E$,
$\frac{{j_{k_p } \left| z \right|r}}{{a_{n_{k_p } } }}\mathop  \to \limits_{p \to \infty } 1$.
There is some $n\in \mathbb{N}$ such that $z=c_n$. We have, for large enough $p$, $f_{j_{k_p } } \left( {\frac{{c_n a_{n_{k_p } } }}{{j_{k_p } }}} \right) = 0$
and
$
\frac{{c_n a_{n_{k_p } } }}{{j_{k_p } }}\mathop  \to \limits_{p \to \infty } rz.
$
So in every neighbourhood of $rz$ there is a zero of $f_{j_{k_{p}}}(z)$ for large enough $p$, and by Lemma 3.4 $rz$ is a non$C_0$-point of $S$ (with respect to which $S$ is irreducible) and it is also a non$C_0$-point of $\left\{ {f_{j_k } } \right\}_{k = 1}^\infty $.
\end{proof}
We deduce
\begin{corollary}
For every sequence $\left\{{f_{j_k}}\right\}$ in $F(f)$, there exists a subsequence $S$ which satisfies one of the following possibilities: \newline
(1) $S$ is a $C_1$-sequence, with unique non$C_0$-point, $z=0$\newline
 or \newline
(2) $S$ is (exactly) a $C_\alpha$-sequence and its set of non$C_0$-points is ${0}\cup r\overline{E}$ for some $r > 0$, with respect to each of these points $S$  is irreducible. \newline
\newline
Possibility (1) in corollary 3.6 indeed occurs, as explained in the following remark.
\end{corollary}
\begin{remark}
Let $L>0$. If we set $j_k=[L\sqrt{a_k a_{k+1}}]$, then $z=0$ is the only non$C_0$-point of ${f_{j_k}}$ in $\mathbb{C}$.
This is because for every $r>0$ we have $a_k<j_k r<a_{k+1}$ for large enough $k$,
$\frac{{j_k r }}{a_{k}} \to \infty $ and  $\frac{a_{k+1}}{{j_k}r} \to \infty$.
So none of the options in Lemma 3.2 can occur.
\end{remark}
\end{proof}
We turn now to the proof (or construction of an appropriate function) of the limit cases in Theorem 1.
\section{Constructing the entire function for the case where $\alpha$ has an i.p and $\nu = \infty$ and for the case where $\alpha$ is a limit ordinal number}
These two cases are treated naturally in the same section, since they are treated similarly. The proof here has similiar features to the proof of the main theorem  in \cite{Generating}. First we discuss the case where $\alpha$ has an i.p and $\nu=\infty$. We prefer here a comprehensive way to the proof for the case $\nu < \infty$. It will save technical calculations and contribute, we hope, to the readability of the proof.
\subsection {Constructing $f_{\alpha,\infty}(z)$ for $\alpha$ having i.p}
First we study the nature of the family $F(f)$ when it is a $Q_\alpha$-normal family of order $\infty$.
One of the following situations must occur (due to Definition 2.34) when we have a family $F$ which is $Q_\alpha$-normal of order $\infty$ (and $\alpha-1$ exists) in $D$.\newline
\textbf{(1)} There exists a sequence $S_\infty \subset F$ having infinitely many non$C_{\alpha  - 1}$-points in $D$, with respect to each of which $S_\infty$ is irreducible. These points cannot have a limit point in $D$, since such a point would be a non$C_{\alpha}$-point, with respect to which $S_\infty$ is irreducible.\newline
\textbf{(2)} For $n$ as large as we like, there exists a sequence $S=S_n$ with at least $n$ finite number of non$C_{\alpha-1}$-points, with respect to each of which $S$ is irreducible.
In the case where the family is $F(f)$ and the domain is $D= \Delta$, option (1) cannot occurs. Indeed, suppose that $S_\infty   = \left\{f_{j_k }\right\}  \subset  F(f)$ has infinitely many non$C_{\alpha-1}$-points $\left\{ {b_m } \right\}_{m = 1}^\infty $ in $\Delta$, with respect to each of which $S_\infty$ is irreducible. Then let $b_0$ be a limit point of these points. According to (1) we must have $\left|b_0\right|=1$, but then $b_0/2\in\Delta$ is a non$C_\alpha$-point of the sequence
$\hat S_\infty  : = \{ f_{2j_k } \}  \subset F(f) $
with respect to which $\hat{S}_\infty$ is irreducible, and it is again a contradiction. So we must rely on option (2) in our construction. We will in fact construct a function $f=f_{\alpha,\infty}$ such that for every $1 \leq n<\infty$ and every $0<r$ corresponds to a sequence in $F(f)$ with exactly $n$ non$C_{\alpha-1}$-points on $\left\{\left| z \right|=r\right\}$, with respect to each of which this sequence in irreducible.
Take $\zeta_0=e^{(i\pi/2)}$ and let $\left\{ {\theta _n } \right\}_{n = 1}^\infty $ be an increasing sequence of positive numbers $\theta_n \nearrow \pi/2$, and set for every $n\geq 1$
$z_n  = e^{i\theta _n } $
($z_n\rightarrow\zeta_0$) and adjust a sector
\begin{equation}
v_n  = S(\theta _n ,\alpha _n )
\end{equation}
 in such a way that $\left\{\Gamma_n:n\geq 1\right\}$ are strongly disjoint, where
\begin{equation}
\Gamma_n:=\overline{v_n\cap\partial\Delta}.
\end{equation}
 Now, by Lemma 2.36, there exists, for every $n \geq 1$, a countable set $E_n=E(\alpha,n,\Gamma_n)$ such that
\begin{equation}
\left( {E_n } \right)_{\partial \Delta }^{(\alpha  - 1)}  = \{ \zeta _1^{(n)} ,\zeta _2^{(n)} ,...,\zeta _n^{(n)} \}.
\end{equation}
We arrange $E_n$ as a sequence $E_n  = \left\{ {z_i^{(n)} } \right\}_{i = 1}^n $. Consider the sequence $\left\{a_n\right\}_{n=1}^\infty$ that satisfies (3.2),(3.3) and (3.4). We introduce now an increasing sequence of positive numbers,$\left\{a_j^{(n)} :n \ge 1\,\,,1 \le j \le n\right\} $ defined as follows: $a_1^{(1)}=a_1;a_1^{(2)}=a_2,a_2^{(2)}=a_3;a_1^{(3)}=a_4,a_2^{(3)}=a_5,a_3^{(3)}=a_6;...$
Define now the following sequence of complex numbers:
\begin{equation}
\begin{array}{*{20}c}
   {a_1^{(1)} z_1^{(1)} ;a_1^{(2)} z_1^{(1)} ,a_1^{(2)} z_2^{(1)} ,a_2^{(2)} z_1^{(2)} ,a_2^{(2)} z_2^{(2)} ;...}  \\
   {a_1^{(n)} z_1^{(1)} ,a_1^{(n)} z_2^{(1)} ,...,a_1^{(n)} z_n^{(1)} ,a_2^{(n)} z_1^{(2)} ,a_2^{(n)} z_2^{(2)} ,...,a_2^{(n)} z_n^{(2)} ,}  \\
   {....a_n^{(n)} z_1^{(n)} ,a_n^{(n)} z_2^{(n)} ,...,a_n^{(n)} z_n^{(n)} ;....}  \\
\end{array}
\end{equation}
The elements of (4.4), by their order defined to be $ \left\{ {b_l } \right\}_{l = 1}^\infty $ and we define similarly to (3.7)
$f(z) = f_{\alpha ,\infty } (z) = \prod\limits_{l = 1}^\infty  {\left( {1 - \frac{z}{{b_l }}} \right)}$.

For each $l$ we have
$\left| {b_l } \right| = a_{t(l)}^{(n(l))} $
, where $n(l)\mathop  \to \limits_{l \to \infty } \infty $ and $1\leq t(l)\leq n(l)$. Also $S(l)=n$ will indicate that $\left|b_l\right|=a_n$. Observe that now on $\left\{ {\left| z \right| = a_n } \right\}$ there are less then $n$ zeros of $f(z)$ for $n \geq 3$. Define
\begin{equation}
E = \bigcup\limits_{n = 1}^\infty  {E_n }.
\end{equation}
 Since ${\Gamma_n:n \geq 1}$ are strongly disjoint, and by the parallel property  guaranteed by Lemma 2.36 we have $E \cap E_{_{\partial D} }^{(1)}  = \emptyset $ and by Lemma 2.12
$
E_{\partial D}^{(\alpha  - 1)}  = \mathop  \cup \limits_{n = 1}^\infty  \{ \zeta _1^{(n)} ,...,\zeta _n^{(n)} \}  \cup \left\{ \zeta_0 \right\}
$
and $E_{\partial\Delta}^{(\alpha)}={\zeta_0}$.
Now, in an analogous way to what we have done in Section 3, we describe the steps, in order to prove that $f(z)$ fulfills the requirements of Theorem~1.
First, Lemma~3.1 is true for our $f_{\alpha,\infty}$. This is because for every fixed $k$
\begin{equation}
\left| {l:n \le S(l) \le n + k\} } \right|\mathop  = \limits_{n \to \infty } o(\left| {\{ l:1 \le S(l) \le n - 1\} } \right|).
\end{equation}
The proof goes exactly the same as  the proof of Lemma~3.1. Lemma~3.2 will also be formulated in the same way, and the proof is very similar.
The only difference is that now the number of zeros of $f(z)$ on $\left\{{\left|z\right|=a_n}\right\}$  is less than $n$ (for $n\geq 3$).
Thus in the calculations in (A),(B) and (C) in the proof of the lemma, instead of $n_k$, there will be another smaller natural number, but by (4.6) we will get the contradiction in the same fashion. The proof of the additional part of this Lemma goes the same. Also in this case we derive already by Lemma 3.2 that $F(f)$ is $Q_\alpha$-normal. Indeed, assume without loss of generality that option (1) in Lemma 3.2 holds, that is
$\frac{{j_{k_p } \left| {\eta _0 } \right|}}{{a_{n_{k_p } } }} \to 1^ +$.
We write $a_{n_{k_p } }  = a_{t\left( {q_p } \right)}^{\left( {n\left( {q_p } \right)} \right)} $. If there is some $t_0$, $1\leq t_0 < \infty$, such that,
\begin{equation}
t(q_p) = t_0
\end{equation}
  for infinitely many p's, say $ \left\{ {p_s } \right\}_{s = 1}^\infty $, then it is easy to see that the \linebreak sequence $S = \left\{ {f_{j_{k_{p_s } } } } \right\}_{s = 1}^\infty $ is a $C_\alpha$-sequence whose set of non$C_0$-points is exactly
$\left| {\eta _0 } \right|\bar E_{t_0 }  \cup \left\{ 0 \right\}$, and with respect to each point in this set $S$ is \linebreak irreducible.
By Lemma 2.20, the non$C_{\alpha-1}$-points of $S$ are \linebreak
$\left| {\eta _0 } \right|\left\{ {\zeta _1^{\left( {t_0 } \right)} ,\zeta _2^{\left( {t_0 } \right)} ,...\zeta _{t_0 }^{\left(t_0 \right) } } \right\}$, all in the sector ${\rm{V}}_{t_0 } $. If, on the other hand, there is no $t_0$ to satisfy (4.7), that is,
$t(q_p )\mathop  \to \limits_{p \to \infty } \infty $, then no $V_n$ contains non$C_0$-points of $\left\{ {f_{j_{k_p } } } \right\}$, and in this case $ \left\{ {f_{j_{k_p } } } \right\} $ is a $C_1$-sequence, with exactly two non$C_0$-points $z_1=0,z_2=\left|\nu_0\right|\cdot\zeta_0$, with respect to each of which $\left\{ {f_{j_{k_p } } } \right\}$ is irreducible, and so $F{(f)}$ is $Q_\alpha$-normal family in $\Delta$.
We observe that although $ E_{\partial \Delta }^{(\alpha )}  = \left\{\zeta _0\right\}  \ne \emptyset $, $F(f)$ is still a $Q_\alpha$-normal family (cf. Lemma~2.20). This is so since for every $n$  the arc $\left\{{\left|z\right|=a_n}\right\}$ contains zeros of $f(z)$ only from one sector $V_t$, $t=t(n)$.\newline
\newline
There is also an analogue to Lemma 3.3 .
\begin{lemma} 
Let $0<r<1$ and $1 \leq t< \infty$. Then there is a sequence in $F(f)$ for which the set on non$C_0$-points is exactly ${0}\cup r\overline{E_t}$.
\end{lemma}
The proof of Lemma 4.1 (to be called Lemma 3.3$^\infty$) is analogous to the proof of Lemma 3.3, but instead of (3.15) we put
\begin {equation}
j_k  = j_k (r,t) = \left[ {\frac{{a_t^{(k)} }}{r} + 1} \right],\quad k \ge t.
\end{equation}
This is because the $(k-t+1)$'th set of zeros in the sector $V_t$ lies on the arc $ \left\{\left| z \right| = a_t ^{\left( k \right)}\right\} $. Observe that as $ k \rightarrow \infty $, the ratio between $a_t^{(k)}$ to the radii of the next arc on which $f(z)$ has zeros tending to zero, but this arc $\left\{\left| z \right| = a_{t + 1}^{(k)}\right\}$ does not contain zeros of $f(z)$ in $V_{t+1}$ and not in $V_t$.
This Lemma together with (4.3) shows that $F(f)$ is not $Q_\alpha$-normal of any finite order and completes the proof of the assertion of Theorem 1 for this case.

Concerning the other assertions about the nature of non$C_0$-sequences in $F(f)$, they all also have analogues in our case.
Lemma 3.4 is true with exactly the same formulation and proof.\newline

The analogue of Lemma 3.5 (to be called $3.5^\infty$) is
\begin{lemma} 
If $rz^{*}$ is a non$C_0$-point of $\left\{ {f_{j_k } } \right\}_{k = 1}^\infty $, where $z^{*}\in\bar{E}_t$, for some $1 \leq t< \infty$, then $rz$ is also a non$C_0$-point of $\left\{ {f_{j_k } } \right\}_{k = 1}^\infty $ for every $z\in\bar{E}_t$.
\end{lemma}
The proof goes along the same lines, with the obvious change of notation. The conclusion is the analogue to Corollary 3.6 (to be called corollary $3.6^\infty$):
\begin{corollary}
For every sequence $\left\{ {f_{j_k } } \right\}$ in $F(f)$, there exists a subsequence $S$, which satisfies one of the following possibilities: \newline
(1) $S$ is a $C_1$-sequence, whith unique non$C_0$-point, $z=0$. or \newline
(2) There is a $1\leq t$  and $r>0$ such that $S$ is a $C_\alpha$-sequence whose set of non$C_0$-points is exactly ${0}\cup r\bar{E}_t$ and with respect to each point of this set $S$ is irreducible.
\end{corollary}
Remark 3.7 is also true, and should be formulated in the same way. Of course, it is related to Corollary 4.3 (corollary $3.6^{\infty}$) and not to Corollary 3.6.\newline

We turn now to the proof of the last case of Theorem 1
\subsection {Constructing $f_\alpha(z)$ when $\alpha$ is a limit ordinal number}
This case is extremely similar to the case where $\alpha-1$ exists and $\nu= \infty$. First, we have a preliminary discussion as we had in section 4.1 on the nature of a $Q_\alpha$-normal family when $\alpha$ is a limit ordinal number. This discussion exposes the similarity between the nature of the family $F(f)$ in this case, to the nature of the family $F(f)$ in the case where $\alpha-1$ exists and $\nu= \infty$. It also suggests the same way for constructing  the corresponding entire function.
Indeed, if $\alpha$ is a limit ordinal number and $F$ is a $Q_\alpha$-normal family in some domain $D$, but not $Q_\beta$-normal there for any $\beta < \alpha$, then one of the following situations must occur according to Definition 2.34:\newline
\textbf{(1)} There exists a sequence $S_\alpha$ in $F$ and an increasing sequence of ordinal numbers $\left\{ {\beta(m)} \right\}_{m = 1}^\infty$, unbounded in $\Gamma(\alpha)$ (i.e., the least upper bound of $\left\{ {\beta(m)} \right\}_{m = 1}^\infty $ is $\alpha$), and to each such $m$ corresponds $z_{\beta(m)} \in D$ which is a non$C_{\beta(m)}$-point of $S_\alpha$, with respect to which $S_\alpha$ is irreducible. The sequence$\left\{ {z_{\beta(m)} } \right\}_{m = 1}^\infty $ cannot have a limit point in $D$, since such a point would be according to Definition 2.34, a non$C_\alpha$-point, with respect to which $S_\alpha$ is irreducible, a contradiction to the $Q_\alpha$-normality of $F$.\newline
\textbf{(2)} For unbounded values of $\beta\in\Gamma(\alpha)$, there exist for each of them a sequence, $S_\beta$ and a point $z_{\beta} \in D$ which is a non$C_\beta$-point of $S_\beta$ with respect to which $S_\beta$ is irreducible.
In the case where the family is $ F = F(f)$ and $D= \Delta$, option (1) cannot occur. Indeed, suppose that $S_\alpha= \left\{f_{j_k }\right\}$ in $F(f)$, and there are infintely many points $\left\{ {z_{\beta(m)}} \right\}_{m = 1}^\infty $ like in (1), where for each $m$, $z_{\beta(m)}\in\Delta$ is a non$C_{\beta(m)}$-point of $S_\alpha$ with respect to which $S_\alpha$ is irreducible and $\left\{ {\beta (m)} \right\}_{m = 1}^\infty $ is unbounded in $\Gamma(\alpha)$. Then let $z_0$ be a limit point of $\left\{z_{\beta(m)}\right\}$. According to (1) we must have $z_0 \in \partial \Delta$. But then $\frac{{z_0 }}{2} \in \Delta $ is a non$C_\alpha$-point of the  sequence $\hat S_\alpha   = f_{2j_k } $ in $F(f)$, with respect to which $\hat{S}_\alpha$ is irreducible and this is again a contradiction. So we must rely on option (2) in our construction.
In fact, we construct a function $f(z)=f_\alpha(z)$ such that for every $\beta \in \Gamma(\alpha)$ and $r > 0$, corresponds a sequence in $F(f)$, with one non$C_\beta$-point on $\left\{{\left| z \right| =r}\right\}$ with respect to which this sequence is irreducible.
To begin this construction, let the enumeration of $\Gamma(\alpha)$ be $\Gamma (\alpha ) = \left\{ {\beta _n } \right\}_{n = 1}^\infty$. $V_n$ and $\Gamma_n$ are defined exactly as in section 4.1 (see (4.1), (4.2)). By Lemma 2.36 for every $n\geq 1$, there is a set $E_n:=E(\beta_n+1,1)\subset\Gamma_n$ such that $\left( {E_n } \right)_\Gamma ^{\left( {\beta _n } \right)}  = \left\{ {z_n } \right\}$.
From this point on, the consruction of $f=f_\alpha$ is exactly the same as the construction of $f=f_{\alpha,\infty}$ in section 4.1, even with the same notations, except that now $f(z)$ is $f_\alpha(z)$ where $ \alpha $ is the limit ordinal number and not $f_{\alpha,\infty}(z)$ as in section 4.1. For every $n \geq 1 $, we use the set $E_n$ in $V_n$ to create $C_{{\beta _n} + 1}$-sequences, which are not $C_{\beta _n }$-sequences and each of them is irreducible with respect to a unique non$C_{\beta_n}$-point, that lies on the ray from the origin through $z_n$. We set $E = \bigcup\limits_{n = 1}^\infty  {E_n } $ (as in (4.5)).

We explain now the various steps in proving that $f_\alpha(z)$ has the desired poperties, doing it in parallel to the previous cases (Section 3 and section 4.1) and in a way that will suffice the reader.
Lemmas 3.1 and 3.2 are formulated the same, with the same proofs, respectively, Lemma 3.2 shows, as in the previous cases, that $F(f)$ is $Q_{\alpha} $-normal in $ \Delta $. The parallel lemma to Lemma~4.1  for this case (to be called Lemma $3.3^{l.o.n})$ is formulated the same and shows that $F(f)$ is not $Q_\beta$-normal for every $\beta < \alpha$. The proof has the same idea, that is, the sequence $\left\{ {f_{j_k } } \right\}$, where $j_k$ is defined by (4.8), is a $C_{\beta(t)+1}$-sequence in $\Delta$ (or in $\mathbb{C}$) and not a $C_{\beta(t)}$-sequence, and it is irreducible with respect to each of its non$C_0$-points that lies on $V_t\cap{\left\{\left|z\right|=r\right\}}$.
Lemma 3.4 is true with the same formulation and proof and so is Lemma 4.2 (that could be called Lemma $3.5^\infty$). The parallel to Corollary 4.3 (or what could be called Corrollary $3.6^\infty$) is
\begin{corollary} 
For every sequence $ \left\{ {f_{j_{k_p } } } \right\}$ in $F(f)$, there exists a subsequence, $S$ which satisfies one of the following possibilities: \newline
\textbf{(1)} $S$ is a $C_1$-sequence, with unique non$C_0$-point, $z=0$ or \newline
\textbf{(2)} There is $t \geq 1$ and $r>0$ such that $\left\{ {f_{j_{k_p } } } \right\}$ is a $C_{\beta(t)+1}$-sequence, whose set of non$C_0$-points is exactly $0 \cup r \bar{E}_t$, with respect to each of which $\left\{ {f_{j_{k_p } } } \right\}$ is irreducible.
\end{corollary}
Also Remark 3.7 has an obvious analogue.\newline
\section {Some remarks}
We introduce now some observations about the entire functions we have created.\newline
\textbf{(1)} In any of the three cases we have treated there is some countable set, $E \subset \partial \Delta$ (where $\overline{E}$ is also countable), such that every non$C_0$-point is lying on a ray through the origin and is a member of $\bar{E}$. We assert that the Julia directions of the corresponding $f(z)$ are also exactly these rays. Morever, in each such Julia direction $f(z)$ attains every complex number infinitely often, with no exceptions. Indeed, if $arg (\eta_0) \notin \left\{{arg z:z \in \bar{E}}\right\}$, then by the proof of Lemma 4.1, $f(z) \Rightarrow \infty$ ($z \rightarrow \infty)$ in $S(arg (\eta_0),\alpha_0/2)$, and thus $L(arg \eta_0)$ is not a Julia direction.
On the other hand, let $z\in\bar{E}$. By Lemma 3.3 (or its analogues), there is a sequence S in $F(f)$, for which  is a non$C_0$-point with respect to which $S$ is irreducible. If there is some $c \in \mathbb{C}$ and $\epsilon >0$ such that $f(z)$ attains $c$ in $S(arg z,\epsilon)$ finitely often, then as in the proof of Lemma 3.4, we can deduce that $z$ is a $C_0$-point of $\left\{ {f_n (z)} \right\}_{n = 1}^\infty $, a contradiction to $z$ being a non$C_0$-point of $S$.\newline
\textbf{(2)} We stated Theorem 1, and gave the corresponding proof and analysis for the family $F(f)=\left\{{f(nz):n \in \mathbb{N}}\right\}$. The same could equally be done to the family $\hat{F}(f) = \left\{ {f(cz)} :c \in \mathbb{C}\right\}$.
Indeed, if $\left\{ {c_k } \right\}_{k = 1}^\infty $ is a sequence of complex numbers $c_k \rightarrow \infty$, then by moving the subsequence we can assume that $\arg c_k \rightarrow \theta_0 \in R$.
Now if we define $n_k=[\left| c_k \right|]$ and $f(z)$ is any meromorphic function on $\mathbb{C}$, then $z_0$ is a $C_0$-point of $\left\{ {f(c_k z)} \right\}_{k = 1}^\infty $ if and only if $z_0$ is a $C_0$-point of $\left\{ {f(n_k e^{i\theta _0 z} )} \right\}$ and the last sequence is just a sequence in the family $F(f)$ acting on the rotated variable $z'=e^{i\theta_0}z$\newline
\textbf{(3)} The order $\rho(f)$ of the entire functions we have created is zero (in all 3 cases). It is very convienient to calculate the exponent of convergence of thes functions. In any of the three cases, the numbers of zeros of $f(z)$ on $\left\{{\left| z\right| =a_n}\right\}$ is less or equal to $n$, thus, for any $\alpha >0$
\begin{equation}
\sum\limits_{l = 1}^\infty  {\frac{1}{{\left| {b_l } \right|^\alpha  }} \le \sum\limits_{n = 1}^\infty  {\frac{n}{{a_n^\alpha  }}} }.
\end{equation}
By conditions  (3.2) and (3.3) we have for large enough $n$, $a_n\geq 2^n$, and thus
\begin{equation}
\sum\limits_{n = 1}^\infty  {\frac{n}{{a_n^\alpha  }}}  < \infty
\end{equation}
and the exponent of convergence is zero and so is $\rho(f)$. \newline
\textbf{(4)} In all three cases the entire function $f$ we have created is in fact $Q_\alpha$-normal in $\mathbb{C}$, not only in $\Delta$. For example, in Lemma 3.1 one can take $\eta_0 \notin c\{L(arg z):z \in \bar{E}$.
Or in Lemma 3.2 one can take $0 \neq \eta_0 \in \mathbb{C}$ and apply exactly the same proofs, respectivly.

\section {An open question}
In \cite{Transfinite} we have shown that given a domain $D$, then for every countable ordinal $\alpha$, there is a family of polynomials which is $Q_\alpha$-normal in $D$, but not $Q_\beta$-normal for any $\beta<\alpha$. Here we have constructed for every countable ordinal $\alpha$ a family of the kind $F(f_\alpha)$ (with suitable $f_\alpha$) having this property.
Let $\alpha_\aleph$ be the first uncountable ordinal (which is of course a limit ordinal). The question naturally arises: does there exist a family which is $Q_{\alpha_\aleph}$-normal (in some domain $D$), but not $Q_\beta$-normal for any countable ordinal number $\beta$?
It seems that the technique we have used so far is not enough to establish the existence of such a family or to negate it. Indeed, our method for constructing such a $Q_\alpha$-normal family, for a countable limit ordinal $\alpha$, is based on creating a countable collection of strongly disjoint sets $\left\{E(\beta)\right\}_{\beta<\alpha}$ in $D$, where for each
$\beta<\alpha$, $E(\beta)_D^{(\beta)}=\emptyset$ and $E(\beta)_D^{(\gamma)}\neq\emptyset$ for every $\gamma<\beta$.
In the uncountable case, the difficulty to use such a technique comes from the fact that there is not enough (topological) `place' in $\mathbb{C}$ for an uncountable collection $\left\{E(\beta)\right\}_{\beta<{\alpha_\aleph}}$ with the property that for every $\beta<\alpha_\aleph$, $E(\beta)_D^{(\beta)}=\emptyset$ and  $E(\beta)_D^{(\gamma)}\neq\emptyset$ for every $\gamma<\beta$.
So the question about the existence of a family which is `exactly' $Q_{\alpha_\aleph}$-normal is open.

\end {document}